\documentclass[draft]{amsart}

\usepackage{ifthen}

\newtheorem{anyprop}{Anyprop}[section]

\newtheorem{theorem}[anyprop]{Theorem}
\newtheorem{lemma}[anyprop]{Lemma}

\newtheorem{corollary}[anyprop]{Corollary}

\theoremstyle{definition}

\newtheorem{example}[anyprop]{Example}

\newtheorem{remark}[anyprop]{Remark}

\newcommand{\NN}{\mathbb{N}}
\newcommand{\ZZ}{\mathbb{Z}}

\newcommand{\RR}{\mathbb{R}}

\newcommand{\PP}{\mathbb{P}}

\newcommand  {\shS}     {\mathcal{S}}

\newcommand  {\Char}    {\operatorname{char}}

\renewcommand{\O}       {\mathcal{O}}

\newcommand  {\Proj}    {\operatorname{Proj}}

\newcommand  {\ra}      {\rightarrow}

\newcommand  {\Syz}     {\operatorname{Syz}}

\newcommand{\komdots}{ , \ldots , }

\newcommand{\lto}{\longrightarrow}

\newcommand  {\td} {\operatorname{td}}
\newcommand{\odd} {\operatorname{odd}}

\theoremstyle{remark}

\numberwithin{equation}{section}

\usepackage{amscd}
\usepackage{amssymb}

\begin{document}
\title[WLP of monomial complete intersections in characteristic $p$]
{A note on the weak Lefschetz property of monomial complete intersections in positive characteristic}

\author[Holger Brenner and Almar Kaid]{Holger Brenner and Almar Kaid}
\address{Universit\"at Osnabr\"uck, Fachbereich 6: Mathematik/Informatik, Albrechtstr. 28a,
49069 Osnabr\"uck, Germany}
\email{hbrenner@uni-osnabrueck.de and akaid@uni-osnabrueck.de}





\begin{abstract}
Let $K$ be an algebraically closed field of characteristic $p > 0$. We apply a theorem of C. Han to give an explicit description
for the weak Lefschetz property of the monomial Artinian complete intersection $A = K[X,Y,Z]/(X^d,Y^d,Z^d)$ in terms of $d$ and $p$.
This answers a question of J. Migliore, R. M. Mir\'{o}-Roig and U. Nagel and, equivalently, characterizes for which characteristics the rank-$2$ syzygy bundle $\Syz(X^d,Y^d,Z^d)$ on $\PP^2$ satisfies the Grauert-M\"ulich theorem. As a corollary we obtain that for $p = 2$ the algebra $A$ 
has the weak Lefschetz property if and only if $d= \left \lfloor \frac{2^t+1}{3} \right \rfloor$ for some positive integer $t$.
This was recently conjectured by J. Li and F. Zanello.
\end{abstract}

\maketitle


\smallskip


\section{Introduction}

Let $R=K[X_0 \komdots X_N]$ be the polynomial ring in $N+1$ variables over an algebraically closed field $K$ and let $f_1 \komdots f_n$ denote
$R_+$-primary homogeneous polynomials in $R$ (i.e., $\sqrt{(f_1 \komdots f_n)}=R_+$). Then the quotient $A:=R/(f_1 \komdots f_n)$ is an Artinian graded $K$-algebra, i.e., $A$ is of the form
$$A = K \oplus A_1 \oplus \ldots \oplus A_s$$ for some integer $s \geq 0$. The algebra $A$ has the \emph{weak Lefschetz property} (abbreviated by WLP) if for every general linear form $\ell \in R_1$ the multiplication maps
$$A_m \stackrel{\cdot \ell} \lto A_{m+1}$$
have maximal rank for $m=0 \komdots s-1$. 

We also associate to the polynomials $f_1 \komdots f_n$ the \emph{syzygy bundle} on $\PP^N=\Proj R$. This vector bundle is given by the short exact sequence
$$0 \lto \Syz(f_1 \komdots f_n) \lto \bigoplus_{i=1}^n \O_{\PP^N}(-d_i) \lto \O_{\PP^N} \lto 0,$$
where $d_i := \deg(f_i)$. If $N = 2$ and $\Char(K) = 0$, we gave in our article \cite{brennerkaidlefschetz} a characterization for the weak Lefschetz property of the Artinian algebra $A$ in terms of the generic splitting type of the syzygy bundle $\Syz(f_1 \komdots f_n)$ (see \cite[Theorem 2.2]{brennerkaidlefschetz}). As a consequence we obtained, using the theorem of Grauert-M\"ulich (see \cite[Theorem 3.0.1] {huybrechtslehn}), the result
of Harima-Migliore-Nagel-Watanabe saying that \emph{every} Artinian complete intersection in $K[X,Y,Z]$ has the weak Lefschetz property (see \cite[Theorem 2.4]{harimamigliore} and \cite[Corollary 2.4]{brennerkaidlefschetz}). The easy examples of the stable syzygy bundles $\shS =\Syz(X^p,Y^p,Z^p)$ over a field of characteristic $p$ show that neither Grauert-M\"ulich ($\shS$ splits on every line $L \subset \PP^2$ as
$\shS|_L \cong \O_L(-p) \oplus \O_L(-2p)$; cf. also the example of L. Ein in \cite[Section 4]{einstable}) nor the theorem of Harima et al. holds in positive characteristic (cf. \cite[Example 7.10]{miglioremironagelmonomial}).

The aim of this paper is to give a numerical characterization of the WLP for monomial Artinian complete intersections $K[X,Y,Z]/(X^d,Y^d,Z^d)$ in positive characteristic. This answers \cite[Question 7.12]{miglioremironagelmonomial} of Migliore-Mir\'{o}-Roig-Nagel and, equivalently, characterizes for which characteristics the rank-$2$ syzygy bundle $\Syz(X^d,Y^d,Z^d)$ on $\PP^2$ satisfies the Grauert-M\"ulich theorem. As a consequence we obtain a proof for the recent conjecture \cite[Conjecture 3.9]{lizanello} of J. Li and F. Zanello.

Besides our geometric approach, the key ingredient for our investigation is a theorem of C. Han which computes the \emph{syzygy gap} for an
ideal of the form $(X^d,Y^d,(X+Y)^d)$ in $K[X,Y]$. 

\section{A characterization of the WLP for monomial Artinian complete intersections in positive characteristic}

The following lemma relates the question of whether the Artinian algebra $A:=K[X,Y,Z]/(X^d,Y^d,Z^d)$ has the weak Lefschetz property to the
generic splitting type of the corresponding syzygy bundle $\Syz(X^d,Y^d,Z^d)$ on the projective plane a la Grauert-M\"ulich.

\begin{lemma}
\label{gmvswlp}
Let $K$ be an algebraically closed field. Then the Artinian complete intersection $K[X,Y,Z]/(X^d,Y^d,Z^d)$ has the weak Lefschetz property
if and only if the syzygy bundle $\shS:=\Syz(X^d,Y^d,Z^d)$ on $\PP^2$ splits on a generic line $L$ as
$\shS|_L \cong \O_L(a) \oplus \O_L(b)$ with $a \geq b$ and $0 \leq a-b \leq 1$.
\end{lemma}

\begin{proof}
Since $\Syz(X^d,Y^d,Z^d)$ is a stable vector bundle (see for instance \cite[Corollary 3.2]{brennerlookingstable}), the assertion of the lemma follows from
\cite[Theorem 2.2]{brennerkaidlefschetz} (note that the proof of this result holds in any characteristic).
\end{proof}

Next, we explain the notion of \emph{syzygy gap} introduced in \cite{monskysyzygygaps} by P. Monsky. Let $K$ be an algebraically closed field and
consider the ideal $I:=(X^{d_1},Y^{d_2},(X+Y)^{d_3})$ in $S:=K[X,Y]$. The minimal graded free resolution of the quotient $S/I$ is given by
$$0 \lto S(a) \oplus S(b) \lto S(-d_1) \oplus S(-d_2) \oplus S(-d_3) \lto S \lto S/I \lto 0,$$
with integers $a,b$, $a \geq b$. The difference $\delta(d_1,d_2,d_3):=a-b$ is called the \emph{syzygy gap} and constitutes a function
$\delta: \NN^3 \ra \NN$. It is easy to see that $a+b = -(d_1+d_2+d_3)$ and hence $\delta(d_1,d_2,d_3) \equiv d_1 + d_2 +d_3 \mod 2$.

\begin{corollary}
\label{relationwlpgmhan}
Let $K$ be an algebraically closed field $($\,of any characteristic\,$)$, $A=K[X,Y,Z]/(X^d,Y^d,Z^d)$ and denote by $\shS = \Syz(X^d,Y^d,Z^d)$ the corresponding syzygy bundle. Then the following conditions are equivalent.
\begin{enumerate}
\item The algebra $A$ has the weak Lefschetz property.
\item The bundle $\shS$ splits on a generic line $L$ as $\shS|_L \cong \O_L(a) \oplus \O_L(b)$ with $a \geq b$ and $0 \leq a-b \leq 1$ $($\,i.e., the theorem of Grauert-M\"ulich holds\,$)$.
\item We have $\delta(d,d,d)\leq 1$.
\end{enumerate}
\end{corollary}

\begin{proof}
The equivalence $(1) \Leftrightarrow (2)$ is proved in Lemma \ref{gmvswlp}.

\noindent $(2) \Leftrightarrow (3)$. If we want to compute the splitting type of $\shS$ on a line $L$ given by the equation $Z=uX+vY$ with coefficients $u,v \in K$, $u,v \neq 0$ (in particular this holds for a generic line), we can  assume without loss of generality that $u=v=1$. Hence computing the generic splitting type of $\shS$ is the same as computing the syzygy gap $\delta(d,d,d)=a-b$.
\end{proof}

We denote by $\delta^*: ~[0, \infty)^3 \ra [0,\infty)$ the continuous continuation of
$\delta$; see \cite[Definition 19 and the following]{monskysyzygygaps} for this function and some of its properties. We set
$$L_{\odd}:=\{(u_1,u_2,u_3) \in \ZZ^3: \sum_{i=1}^3 u_i \mbox{
odd}\} \subset \ZZ^3.$$
An element $u=(u_1,u_2,u_3) \in \ZZ^3$ belongs to $L_{\odd}$ if and only if all entries of $u$ are odd or if there is only
one odd entry $u_i$, $i \in \{1,2,3\}$. Further, we denote by $\td$ the \emph{taxi-cab
distance} in $\RR^3$ defined as $\td(v,w):= \sum_{i=1}^3|v_i-w_i|$
for triples $v=(v_1,v_2,v_3), w=(w_1,w_2,w_3) \in \RR^3$.

The following theorem due to C. Han yields an effective way to compute $\delta^*$ for a
given triple $v=(v_1,v_2,v_3) \in [0, \infty)^3$.

\begin{theorem}[Han]
\label{hanstheorem}
Let $K$ be an algebraically closed field of characteristic $p>0$ and assume the entries of
$v=(v_1,v_2,v_3) \in [0, \infty)^3$ satisfy $v_1 \leq v_2 \leq v_3$ and $v_3 < v_1 + v_2$.
If there exists $s \in \ZZ$ and a triple $u=(u_1,u_2,u_3) \in
L_{\odd}$ such that $m:=\td(p^s v,u)<1$, then there exists such a pair
$s,u$ with minimal $s$. With these data $s,u$ and $m$ we have
$$\delta^*(v)=p^{-s}(1-m).$$ If no such pair exists, then
$\delta^*(v)=0$.
\end{theorem}

\begin{proof}
See \cite[Theorems 2.25 and 2.29]{hanthesis} or \cite[Corollary
23]{monskysyzygygaps} for an easier proof.
\end{proof}

\begin{lemma}
\label{distancelemma}
Let $d \in \NN_+$ and $p$ be a prime number. Then the following conditions are equivalent.
\begin{enumerate}
\item There exists $k \in \NN$ and $n \in \NN$ such that
\[  \frac{3d}{6k+2} >  p^n > \frac{3d}{6k+4} \, .\]
\item
There exists an odd number $u \in \NN$ and $s \in \ZZ,\, s \leq 0$, such that
\[ u- \frac{1}{3} < d p^s < u + \frac{1}{3} \, .\]
\item
There exists an integer $s$, $s \leq 0$, such that the taxi-cab distance of $(dp^{s}, dp^{s} ,dp^{s})$ to some point in $L_{\odd}$ is $<1$.
\end{enumerate}
\end{lemma}

\begin{proof}
To proof the equivalence between (1) and (2) we set $s=-n$ and $u=2k+1$. The condition in (1) is equivalent with
\[ \frac{3}{3u-1} > \frac{p^n}{d} > \frac{3}{3u+1}  \]
and by inverting it is equivalent with
\[  u - \frac{1}{3} < d {p^s} < u+ \frac{1}{3}  \, .\]

If (2) is true, then we have $(u,u,u) \in L_{\rm odd}$ and the taxi-cab distance between $(dp^{s}, dp^{s} ,dp^{s})$ and $(u,u,u)$ is $<1$. On the other hand, the distance of a point on the diagonal to any point in $L_{\odd}$ outside the diagonal is at least $1$, so we only have to consider points on the diagonal.
\end{proof}

\begin{lemma}
\label{maximalnlemma}
Let $d \in \NN_+$ and $p$ be a prime number. Suppose that there exists $0 \leq n' < n$ and $k',k \in \NN$ such that
\[  \frac{3d-1}{6k'+2}  > p^{n'} > \frac{3d+1}{6k'+4} \]
and
\[  \frac{3d}{6k+2}  > p^{n} > \frac{3d}{6k+4}  \, .  \]
Then
\[  \frac{3d-1}{6k+2}  > p^{n} > \frac{3d+1}{6k+4}  \, .  \]
\end{lemma}

\begin{proof}
Otherwise we would have either
\[  \frac{3d+1}{6k+4}  \geq p^{n} > \frac{3d}{6k+4}  \]
or
\[  \frac{3d}{6k+2}  > p^{n} \geq \frac{3d-1}{6k+ 2}  \, .\]
This gives either
\[p^{n} (6k+4) =3d+1\]
or
\[p^{n} (6k+2) =3d-1 \, .\]
We plug this in the first inequality and get in the first case
\[ \frac{ p^{n} (6k+4) -2  }{6k'+2}  > p^{n'} > \frac{p^{n} (6k+4)   }{6k'+4}   \]
and by dividing through $p^{n}$ we get
\[ \frac{3k+2   - \frac{1}{p^{n} }}{3k'+1} >  p^{n'-n} > \frac{3k+2}{3k'+2} \, . \]
By inverting we obtain
\[ \frac{3k'+1}{3k+2   - \frac{1}{p^{n} }} <  p^{n-n'} < \frac{3k'+2}{3k+2} \, . \]
From the right hand side we get $  p^{n-n'} \leq  \frac{3k'+1}{3k+2} $ which yields the contradiction
\[ \frac{3k'+1}{3k+2   - \frac{1}{p^{n} }} < \frac{3k'+1}{3k+2} \, . \]

In the second case we obtain
\[  \frac{p^{n} (6k+2)}{6k'+2}  > p^{n'} > \frac{p^{n} (6k+2) +2 }{6k'+4} \]
and similar manipulations yield a contradiction.
\end{proof}

The following theorem gives an explicit answer to \cite[Question 7.12]{miglioremironagelmonomial}. This question was also answered in \cite[Corollary 3.6]{lizanello} but in a less explicit way.

\begin{theorem}
\label{maintheorem}
Let $K$ be a field of characteristic $p>0$ and consider the monomial Artinian complete intersection $A:=K[X,Y,Z]/(X^d,Y^d,Z^d)$.
Then the following holds:
\begin{enumerate}
\item If $d$ is even, then $A$ does not have the weak Lefschetz property if and only if there exists a $k \in \NN$ and an $n \in \NN_+$ such that
$$\frac{3d}{6k+2} > p^n > \frac{3d}{6k+4}.$$
\item If $d$ is odd, then $A$ does not have the weak Lefschetz property if and only if there exists a $k \in \NN$ and an $n \in \NN_+$ such that
$$\frac{3d-1}{6k+2} > p^n > \frac{3d+1}{6k+4}.$$
\end{enumerate}
\end{theorem}

\begin{proof}
We prove $(1)$. Assume that we have
$$\frac{3d}{6k+2} > p^n > \frac{3d}{6k+4}$$
for some $k \in \NN$ and $n \in \NN_+$. We set $s:=-n$, $u:=2k+1$. Then we have $m:=\td(p^{s}(d,d,d),(u,u,u))<1$ by Lemma \ref{distancelemma} and
hence $\delta^*(d,d,d)=p^{-s}(1-m)>0$. Since $\delta^*(d,d,d)=a-b$ and $a+b =-3d$ we must have $\delta^*(d,d,d) \geq 2$.
We apply Corollary \ref{relationwlpgmhan} and see that $A$ does not have the WLP.

\smallskip

Now we assume that the numerical condition does not hold. Then by Lemma \ref{distancelemma} there is no $s \leq 0$ such that
the taxi-cab distance from $p^s(d,d,d)$ to an element $(u,u,u) \in L_{\odd}$ is $<1$. This is also true for $s>0$ since $d$ is even. Hence it follows from Han's Theorem \ref{hanstheorem} that $\delta^*(d,d,d) = 0$ which implies by Corollary \ref{relationwlpgmhan} the WLP for the algebra $A$.

\smallskip

Next we prove $(2)$. First we remark that, since $d$ is odd, the condition
$$\frac{3d}{6k+2}  > p^{n} > \frac{3d}{6k+4}$$
is always fulfilled for $n=0$ and $k$ such that $d=2k+1$. We choose $n>0$ maximal such that
$$\frac{3d}{6k+2}  > p^{n} > \frac{3d}{6k+4}$$ holds for some $k$. Hence we can apply Han's Theorem \ref{hanstheorem} with $s:=-n$ (minimal) and $u:=2k+1$ to compute the syzygy gap.

\smallskip

Suppose that the numerical condition of part $(2)$ is fulfilled for some $k' \in \NN$ and $n' \in \NN_+$.  According to Lemma \ref{maximalnlemma} we may assume that this condition also holds for the chosen (maximal) $n$, hence 
$$\frac{3d-1}{6k+2} > p^n > \frac{3d+1}{6k+4}\,.$$
Then we have in particular
$$u-\frac{1}{3} = \frac{6k+2}{3} < dp^s < \frac{6k+4}{3} = u+ \frac{1}{3}$$
by Lemma \ref{distancelemma}. Now we distinguish two cases.

\smallskip

\noindent \textbf{Case 1:} Let $u > dp^s$. Then the taxi-cab distance from $p^s(d,d,d)$ to the element $(u,u,u) \in L_{\odd}$ equals
$$m:=\td(p^s(d,d,d),(u,u,u)) =3(u-dp^s)$$
and we have $m < 1$ (by Lemma \ref{distancelemma}). So we obtain for the syzygy gap:
\begin{eqnarray*}
\delta^*(d,d,d)&=& p^{-s}(1 - m)\\ &=& p^{-s} (1-3u+3dp^s)\\&=& p^{-s}(1-3u)+3d\\&=&-p^n (6k+2) +3d\\&>& -(3d-1) +3d\\&=&1.
\end{eqnarray*}
Therefore the syzygy gap is indeed $\geq 3$. Hence it follows from Corollary \ref{relationwlpgmhan} that $A$ does not have the WLP.

\smallskip

\noindent \textbf{Case 2:} Let $u \leq  dp^s$. Then we obtain
$$m:=\td(p^s(d,d,d),(u,u,u)) =3(dp^s-u)$$
which is again $<1$. So we can estimate the syzygy gap as follows:
\begin{eqnarray*}
\delta^*(d,d,d)&=& p^{-s}(1 - m)\\&=& p^{-s}(1 + 3u-3dp^s)\\&=& (1+3u)p^{-s}-3d\\&=&(6k+4)p^n-3d\\&>&3d+1 -3d\\&=&1.
\end{eqnarray*}
Again we conclude that $A$ does not have the WLP.

\smallskip

Next suppose that the numerical condition of part $(2)$ does not hold. Then we have either 
$$\frac{3d+1}{6k+4} \geq p^n > \frac{3d}{6k+4} \mbox { or } \frac{3d}{6k+2} > p^n \geq \frac{3d-1}{6k+2},$$
where $n$ and $k$ are chosen as in the beginning of the proof of part $(2)$.

\smallskip

\noindent \textbf{Case 1:} Let $\frac{3d+1}{6k+4} \geq p^n > \frac{3d}{6k+4}$. Then we even have $$p^n(6k+4)=p^n(3u+1)=3d+1.$$
Since $$\frac{d}{u} =\frac{3d}{3u} > \frac{3d+1}{3u+1}= \frac{3d+1}{6k+4} = p^n,$$
we have $dp^s > u$. So we obtain
$$m:=\td(p^s(d,d,d),(u,u,u))= 3(dp^s-u)$$
which is $<1$. This gives:
\begin{eqnarray*}
\delta^*(d,d,d)&=& p^{-s}(1 - m)\\&=& p^{-s}(1 + 3u -3dp^s)\\&=& p^{-s}(1+3u) -3d\\&=& p^n(3u+1) - 3d\\&=& 1.
\end{eqnarray*}
Hence $A$ has the WLP by Corollary \ref{relationwlpgmhan}.

\smallskip

\noindent \textbf{Case 2:} Let $\frac{3d}{6k+2} > p^n \geq \frac{3d-1}{6k+2}$. This implies $$p^n(6k+2)=-p^n(1-3u)=3d-1.$$ Since
$$p^n =\frac{3d-1}{6k+2} = \frac{3d-1}{3u-1} < \frac{d}{u},$$
we now have $u > dp^s$. Hence
$$m:=\td(p^s(d,d,d),(u,u,u))= 3(u-dp^s)<1.$$
Once again we get
\begin{eqnarray*}
\delta^*(d,d,d)&=& p^{-s}(1 - m)\\&=&p^{-s}(1-3u-3dp^s)\\&=&p^{n}(1-3u)+3d\\&=& -(3d-1)+3d\\&=& 1.
\end{eqnarray*}
We conclude as above that $A$ has the WLP.
\end{proof}

As a corollary we obtain \cite[Conjecture 3.9]{lizanello}.

\begin{corollary}
\label{zanelloconjecture}
Let $K$ be a field of characteristic $2$. Then the Artinian complete intersection $A:=K[X,Y,Z]/(X^d,Y^d,Z^d)$ has the weak Lefschetz property if and only if $d= \left \lfloor \frac{2^t+1}{3} \right \rfloor$ for some positive integer $t$.
\end{corollary}

\begin{proof}
Let $n \in \NN$ such that
$$\frac{3d}{2} > 2^n > \frac{3d}{4}$$
(note that there is only one such $n$ since $\frac{3d}{4}$ is the half of $\frac{3d}{2}$).
This $n$ corres\-ponds to $k=0$ and is the exponent we have to consider by Theorem \ref{maintheorem}. So it follows from part $(1)$ of Theorem \ref{maintheorem} that the algebra $A$ never enjoys the WLP for $d$ even. So we may
assume that $d$ is odd. If
$$\frac{3d-1}{2} > 2^n > \frac{3d+1}{4}$$
holds then again $A$ does not have the WLP. So $A$ does have the WLP if either
$$\frac{3d-1}{2}\leq 2^n < \frac{3d}{2} \mbox{ or }  \frac{3d}{2} < 2^{n+1} \leq \frac{3d+1}{2}$$
holds, i.e., if we have either $3d-1 = 2^{n+1}$ or $3d+1 = 2^{n+2}$. This gives the assertion of the corollary.
\end{proof}

\begin{remark}
\label{evendchar2}
As remarked in \cite{lizanello} and indicated in our proof, Corollary \ref{zanelloconjecture} implies that the monomial complete intersection $K[X,Y,Z]/(X^d,Y^d,Z^d)$ does \emph{not} have the WLP in characteristic $2$ if $d$ is even.
\end{remark}

Theorem \ref{maintheorem} implies in particular that for given $d$ the weak Lefschetz property might only fail in characteristic $p \leq \frac{3}{2} d$. It is easy to generate the list of exceptional characteristics with the help of this numerical criterion.

\begin{corollary}
\label{p teilt d}
Let $d$ be odd and let $p$ be a prime factor of $d$. Then the Artinian algebra $K[X,Y,Z]/(X^d,Y^d,Z^d)$ does not have the weak Lefschetz property in characteristic $p$.
\end{corollary}

\begin{proof}
We write $d=p^n u$ with $u=2k+1$ odd, $n \geq 1$.
Then
\[ p^n = \frac{d}{2k+1} = \frac{3d}{6k+3} \, . \]
Since the numerator is larger than the denominator, this number is strictly between $\frac{3d+1}{6k+3+1}$ and $\frac{3d-1}{6k+3-1}$, so this fulfills the condition of Theorem \ref{maintheorem}(2).
\end{proof}

Remark \ref{evendchar2} and Corollary \ref{p teilt d} imply that only for $d=1$ the WLP holds in all characteristics. We will see in the examples below that for $d$ even the weak Lefschetz property can hold in characteristics dividing $d$ (but not in characteristic $2$).

\begin{example}
We consider $d$ even and determine the exceptional prime numbers (here we mean by \emph{exceptional} that the Artinian complete intersection $A=K[X,Y,Z]/(X^d,Y^d,Z^d)$ does not enjoy the WLP in these characteristics).

$d=2$. The only exceptional prime number is $2$.

$d=4$. The condition for $k=0$ is $ 12/2=6 > p^n > 12/4=3$, hence the exceptional prime numbers are $2$ and $5$ (no larger $k$ have to be considered).

$d=6$. For $k=0$ we get $9 > p^n > 4.5 $, which yields the exceptional primes $2,5,7$ (no larger $k$). The prime number $3$ divides $d$, but the weak Lefschetz property does hold in characteristic $3$.

$d=8$. For $k=0$ we get $12 > p^n > 6$, which yields the exceptional primes $2,3,7,11$ (no larger $k$).

$d=10$. For $k=0$ we get the exceptional primes $2,3,11,13$ (no larger $k$).

$d=12$. For $k=0$ we get the exceptional primes $2,11,13,17$ (no larger $k$).

$d=14$. For $k=0$ we get the condition $21 > p^n > 10.5 $, which yields the exceptional primes $2,11,13,17,19$. For $k=1$ we get the condition $\frac{42}{8} > p^n > \frac{42}{10}$, which yields $p=5$.

$d=20$. For $k=0$ we get the exceptional primes $2,3,5,17,19,23,29$ and for $k=1$ we also get $7$. Note that $5$ does divide $d$ and the algebra does not have the weak Lefschetz property in characteristic $5$.
\end{example}

\begin{example}
We consider $d$ odd and determine the exceptional prime numbers.

$d=1$. For $k=0$ we get the condition $ 1> p^n >1$, which has no solution, hence $K[X,Y,Z]/(X,Y,Z) \cong K$ has the weak Lefschetz property in every characteristic, which is clear anyway.

$d=3$. The condition for $k=0$ is $ 8/2= 4 > p^n > 10/4=2.5$, hence the only exceptional prime number is $3$ (no larger $k$ have to be considered).

$d=5$. For $k=0$ we get $7 > p^n > 4 $, which yields the only exceptional prime $5$ (no larger $k$). The prime number $7$ fulfills $7 = \frac{14}{2}= \frac{3d-1}{2}$, which corresponds to the second case in the proof of Lemma \ref{maximalnlemma}. For $p=7$ the Han number is $s=-1$, but the syzygy gap is $1$ and not $3$.

$d=7$. For $k=0$ we get $10 > p^n > 5.5$, which yields the exceptional primes $2,3,7$ (no larger $k$).

$d=9$. For $k=0$ we get the exceptional primes $2,3$ and $11$ (no larger $k$).

$d=31$. For $k=0$ we get the condition $46 > p^n > 23.5$, which yields the exceptional primes $2,3,5,29,31,37,41,43$. For $k=1$
we get the condition $ \frac{92}{8} > p^n > \frac{94}{10}$, which yields also $p=11$.
\end{example}

\bibliographystyle{amsplain}

\end{document}